\documentclass{article}

\usepackage{amssymb,amsfonts,amsmath,latexsym,version, bigints, comment}
\usepackage[latin1]{inputenc}
\usepackage{url}
\newtheorem{thm}{Theorem}[section]
\newtheorem{lem}[thm]{Lemma}
\newtheorem{defn}[thm]{Definition}
\newtheorem{prop}[thm]{Proposition}
\newtheorem{rem}[thm]{Remark}
\newtheorem{rems}[thm]{Remarks}
\newtheorem{ex}[thm]{Example}
\newtheorem{cor}[thm]{Corollary}
\newtheorem{met}[thm]{Method}

\newcommand\Z{{\mathbb Z}}
\newcommand\N{{\mathbb N}}
\newcommand\R{{\mathbb R}}

\newcommand\E{{\mathbb E}}
\newcommand\eps{\varepsilon}
\newcommand{\Cov}{\mbox{\rm Cov}}
\newcommand{\Var}{\mbox{\rm Var}}

\author{\small{Julie Fournier, MAP5, universit\'e Paris Descartes and universit\'e Pierre et Marie Curie}}
\date{\small{March 2017}}
\begin{document}

\title{Identification and isotropy characterization of deformed random fields through excursion sets} 

\maketitle

\author 
\date
\paragraph{Abstract}
A deterministic application $\theta\,:\,\mathbb{R}^2\rightarrow\mathbb{R}^2$ deforms bijectively and regularly the plane and allows to build a deformed random field $X\circ\theta\,:\,\mathbb{R}^2\rightarrow\mathbb{R}$ from a regular, stationary and isotropic random field $X\,:\,\mathbb{R}^2\rightarrow\mathbb{R}$. The deformed field $X\circ\theta$ is in general not isotropic, however we give an explicit characterization of the deformations $\theta$ that preserve the isotropy.
Further assuming that $X$ is Gaussian, we introduce a weak form of isotropy of the field $X\circ\theta$, defined by an invariance property of the mean Euler characteristic of some of its excursion sets. Deformed fields satisfying this property are proved to be strictly isotropic. Besides, assuming that the mean Euler characteristic of excursions sets of $X\circ\theta$ over some basic domains is known, we are able to identify $\theta$.\\
\vspace{1cm}

Deformed fields are a class of non-stationary and non-isotropic fields obtained by deforming a fixed stationary and isotropic random field thanks to a deterministic function that transforms bijectively the index set.
Deformed fields respond to the need to model spatial and physical phenomena that are in numerous cases not stationary nor isotropic. To give but one example, they are currently widely used in cosmology to model the cosmic microwave background (CMB) deformed anisotropically by the gravitational lensing effect, with mass reconstruction as an objective \cite{HO02}.

Our framework is two-dimensional: we set $X\,:\,\R^2\rightarrow\R$ the underlying stationary and isotropic field, $\theta\,:\,\R^2\rightarrow\R^2$ a bijective deterministic function and $X_{\theta}=X\circ\theta$ the deformed field. In fact, most studies on the deformed field model deal with dimension two. The reason for this is that it is the simplest case of multi-dimensionality, the results can be illustrated easily thanks to simulations and it still covers a lot of possible applications, particularly in image analysis. For instance, deformed fields are involved in the "shape from texture" issue, that is, the problem of recovering a 3-dimensional textured surface thanks to a 2-dimensional projection \cite{CM02}.   

The model of deformed fields was introduced in 1992 in a spatial statistics framework by Sampson and Guttorp in \cite{SG92}, with only a stationarity assumption on $X$. It is also studied through the covariance function in \cite{PM99} and in \cite{PS00}. In \cite{ASP15}, the authors investigate the case of a linear deformation with a matrix representation as the product of a diagonal and a rotation matrix, which produces what they call "geometric anisotropy". In \cite{CM03}, the deformed field model is studied as a particular case of a model of deterministic deformation operator applied to a stationary field $X$. A lot of papers also propose methods to estimate $\theta$, as we will see a little further on in this introduction, when we come to our own contribution to the estimation matter.

Unless otherwise specified, the kind of stationarity and isotropy that we consider consists in an invariance of the field's law under translations or, respectively, rotations. Even though the underlying field $X$ is stationary and isotropic, a lot of deformations transform the index space $\R^2$ in such a way that the stationarity and/or the isotropy are lost when it comes to the deformed field.
The deformations preserving stationarity are the linear deformations. Concerning isotropy, a natural question arises: which are the deformations $\theta$ that preserve isotropy, for any underlying field $X$ ? It is solved in Section \ref{sectionstatiso}. We give an explicit form for this kind of deformations and we call them spiral deformations. Let us point out here that the question of preserving the isotropy for one fixed underlying field $X$ is different, and it is solved in Section \ref{section-stariso}.

For the rest of the paper, we have in mind the following practical problem: the covariance function of the underlying field $X$ and the deformation $\theta$ are unknown. We try to study and even to identify $\theta$ through observations of some excursion sets of $X_{\theta}$ above fixed levels. 
For this, we add some assumptions on $X$ (Gaussianity, $\mathcal{C}^2$-regularity, non-degeneracy assumptions) and on $\theta$ ($\mathcal{C}^2$-regularity), which are precisely described and justified in Section \ref{sectionnota}, and we focus on the mean Euler characteristic of the excursion sets. The Euler characteristic is an additive topological functional that is defined on a large class of compact sets. Heuristically, the Euler characteristic of a set is determined by its topology: for a two-dimensional compact set, it is the number of connected components minus the number of holes in this set; for a one-dimensional set, it is simply the number of closed intervals that compose the set. Note that a modified version of the Euler characteristic of excursion sets will be more adequate to address our problem. The formulas of the expectation of the (modified) Euler characteristic of an excursion set of $X_{\theta}$ can be found in Section \ref{section-Euler}. 

More precisely, let $T$ be a compact regular set in $\R^2$. We are interested in the Euler characteristic $\chi$ of the excursion set of $X_{\theta}$ restricted to $T$ above a level $u\in\R$, $A_u(X_{\theta},T)=\{t\in T\;/\; X(\theta(t))\geq u\}$. However, we may study equivalently the stationary and isotropic field $X$ on the transformed set $\theta(T)$ or the deformed field $X_{\theta}$ on the set $T$, since
\[\chi(A_u(X_{\theta},T))=\chi(A_u(X,\theta(T))).
\]
The sets $T$ that we consider are basic domains: segments and rectangles in $\R^2$. 

In Section \ref{section-stariso}, we introduce the notion of $\chi$-isotropic deformation: it applies to a deformation $\theta$ such that, for any level $u$ and for any rectangle $T$, $\E[\chi(A_u(X_{\theta},T))]$ does not vary under any rotation of $T$. This is in particular true if the deformed field is isotropic, hence this property can be viewed as a weak notion of isotropy. However, it occurs that this weak notion implies the strong one (isotropy in law), that is, the $\chi$-isotropic deformations are exactly the spiral deformations.

In Section \ref{section-theta-unknown}, we tackle the problem of identifying $\theta$, assuming that we only have at our disposal the mean Euler characteristic of some excursion sets of the deformed field. 
The problem of the estimation of a deformation $\theta$ thanks to the observation of the deformed random field $X_{\theta}$ is originally a spatial statistics problem and it has been studied under different angles since it was introduced.                         
At first, \cite{SG92} used several observations on a sparse grid to estimate $\theta$.
Another approach is to use only one observation of the deformed field on a dense grid; it is adopted in \cite{GP00}, \cite{CM03}, \cite{AS08}, \cite{AC09}, and \cite{FDR15} with an underlying field that is stationary and/or isotropic. These different papers involve convergence results on quadratic variations and quasi-conformal theory. The study in \cite{AG16} applies in particular to deformed fields of the form $\{X(x+\nabla  \eta(x)),\:x\in\R^2\}$, for which an estimation of the deterministic $\eta:\R^2\rightarrow \R$ is proposed; this very model is indeed used in cosmology for the estimation of the gravitational lensing of the CMB \cite{HO02}.

Our approach differs from the previous ones, since our observations are limited to realizations of $X_{\theta}$ over a fixed level, and not to the whole realizations. Our method is closer to the one in \cite{C87}, where the inference of the deformation is based on the size and shape of the deformed field's level curves; however, the author restricts the deformations to linear ones given by symmetric, positive and definite matrices. With our sparse observations, we manage as well as in \cite{AS08} to compute the complex dilatation of $\theta$  up to a conformal map, at every point of the domain. The complex dilatation provides a characterization of the deformation. 

In this paper, we prove four main results. Theorem \ref{equivspiral} states that the deformations preserving isotropy are exactly the spiral deformations. In Theorem \ref{thm_chiiso}, the class of deformations satisfying the invariance condition of the mean Euler characteristic of excursion sets is identified with the spiral deformations. A consequence of this theorem is Corollary \ref{Coro}. Roughly speaking, it states that three notions of preservation of isotropy coincide and correspond to the set of spiral deformations. In Section \ref{section-theta-unknown}, we show how to almost entirely identify $\theta$ through the mean Euler characteristic of its excursion sets over basic domains. The general case is described by Method \ref{propgeneralcase}. To end with, in Section \ref{subsec_spiral}, limiting ourselves to spiral deformations, we finally propose an estimation method based on one single observation of the deformed field.
 \section{Notations and assumptions}
 \label{sectionnota}
 For any compact $A$ in $\R^2$, we write $\dim(A)$ its Hausdorff dimension; if $\dim(A)=1$, we write $|A|_1$ its one-dimensional Hausdorff measure; if $\dim(A)=2$, we write $|A|_2$ its two-dimensional Hausdorff measure and $\partial A$ its frontier.
 
We work in a fixed orthonormal basis in $\R^2$ and we will use the same notation for a linear application defined on $\R^2$ and taking value in $\R^2$ and for its matrix in this basis.  We denote by $O(2)$ (respectively $SO(2)$) the group of orthogonal transformations in $\R^2$ (respectively the group of rotations in $\R^2$), by $\mathbb{T}$ the one-dimensional torus and for any $\alpha\in\mathbb{T}$, by $\rho_{\alpha}$ the rotation of angle $\alpha$.

 For any real $s$, we write $[0,s]=\{x\in\R,\, 0\le x\le s\}$ if $s\geq0$ and $[0,s]=\{x\in\R,\, s\le x\le 0\}$ if $s<0$. We say a set of points $T$ in $\R^2$ is a segment if there exists $(a,b)\in(\R^2)^2$ with $a\neq b$ such that $T=\{a+t(b-a),\,t\in[0,1]\}$. For any $(s,t)\in\R^2$, we write $T(s,t)=[0,s]\times[0,t]$ and we say a set of points $T$ in $\R^2$ is a rectangle if there exist $(s,t)\in\R^2$, $\rho\in SO(2)$ and a translation $\tau$ such that $T=\rho\circ\tau(T(s,t))$.
 
 If $f=(f_1,f_2):\R^2\rightarrow \R^2$, with $f_i:\R^2\rightarrow\R$ for $i\in\{1,2\}$, is a differentiable function, for any $x=(s,t)\in \R^2$, we use the notations $J_f^1(x)$ for the vector $\partial_s f(x)=(\partial_s f_1(x),\partial_s f_2(x))$, $J_f^2(x)$ for the vector $\partial_t f(x)$ and $J_f(x)$ for the Jacobian matrix of $f$ at point $x$. More generally, if $M$ is a $2\times 2$ matrix, for $i\in\{1,2\}$, we write $M^i$ the $i^{\text{th}}$ column of $M$.

Let $X$ be a Gaussian, stationary and isotropic random field, defined on $\R^2$ and taking real values; we write $C \, : \, \R^2 \rightarrow \R$ its covariance function. Since $X$ is stationary, we may assume it is centered. We shall also assume that $C(0)=1$ since if not, we consider the field $\frac{1}{\sqrt{C(0)}}X$ instead of $X$. As for the regularity of $X$, we make the assumption that almost every realization of $X$ is of class $\mathcal{C}^2$ on $\R^2$.
As a consequence, $C$ is of class $\mathcal{C}^4$. We denote by $X'(t)$ (respectively by $X''(t)$) the gradient vector (respectively the Hessian matrix) of $X$ at point $t$ and by $C''(t)$ the Hessian matrix of $C$ at point $t$.   
We assume that for any $t\in\R^2$, the joint distribution of $(X_i'(t),X_{i,j}''(t))_{(i,j)\in\{1,2\}^2,i\le j}$ is not degenerate, in order to be able to apply the mean Euler characteristic of excursion sets formula. Therefore, the covariance matrix of $X'(0)$ is not degenerate; since $X$ is isotropic, there exists $\lambda>0$ such that $\Cov(X'(0))=\lambda\,I_2$. If $\lambda\neq1$, $X_{\theta}$ is nevertheless equal to $\tilde{X}_{\tilde{\theta}}$, with $\tilde{\theta}=\sqrt{\lambda}\theta$ and with $\tilde{X}(\cdot)=X(\sqrt{\lambda}^{-1}\cdot)$ satisfying $\Cov(\tilde{X}'(0))=I_2$. Consequently, without loss of generality, we shall assume that $C''(0)=-I_2$. 

We gather all the assumptions on $X$ that will be in force in Sections \ref{section-Euler}, \ref{section-stariso} and \ref{section-theta-unknown} under the name \textbf{(H)}:
\begin{equation*}
\left\{
\begin{aligned}
&X\text{ is Gaussian,}\\
&X\text{ is stationary and isotropic,}\\
&X\text{ is almost surely of class }\mathcal{C}^2,\\
&\forall t\in\R^2,\text{ the joint distribution of }(X_i'(t),X_{i,j}''(t))_{\underset{i\le j}{(i,j)\in\{1,2\}^2}}\text{ is not degenerate,}\\
&X\text{ is centered, }C(0)=1\text{ and }C''(0)=-I_2.
\end{aligned}
\right.
\end{equation*}

Our ambition in Section \ref{section-theta-unknown} is to identify the deformation $\theta$ assuming that we only have at our disposal the expectation of $\chi(A_u(X_{\theta}, T))$ for different sets $T$ and for a fixed level $u$. However, it is not possible to distinguish between $\theta$ and another deformation $\tilde{\theta}$ such that the random fields $X_{\theta}$ and $X_{\tilde{\theta}}$ have the same law, which defines an equivalence relation between $\theta$ and $\tilde{\theta}$. Because of the stationarity and the isotropy of $X$, the equivalence class of $\theta$ is $\{\tilde{\theta}=\rho\circ\theta+u,\;\rho\in O(2),\;u\in\R^2\}$ and we can only hope to determine one of its representative and not $\theta$ itself. 

Consequently, without loss of generality, we can make the assumption that $\theta(0)=0$. 
If $\theta$ is differentiable, we shall also assume that for any $x\in\R^2$, $\det( J_{\theta}(x))$ is positive or, in other words, that $\theta$ preserves orientation. Indeed, function $x\mapsto \det(J_{\theta}(x))$ is continuous on $\R^2$ and does not take zero value, hence it takes either only positive or only negative values. If for all $x\in\R^2$, $\det(J_{\theta}(x))<0$, we can replace $\theta$ by $\sigma\circ \theta$, where $\sigma\in O(2)$ is the symmetry with respect to the axis of abscissa; then for any $x\in\R^2$, $J_{\sigma\circ\theta}(x)=\sigma\circ J_{\theta}(x)$ and so the Jacobian determinant of $\sigma\circ\theta$ is positive on $\R^2$. 
Those two transformations on $\theta$ (translation along vector $\theta(0)$ and left composition with $\sigma$) only mean that we consider another representative in the equivalence class of $\theta$.
Note that the class of linear as well as tensorial deformations considered as examples in Section \ref{section-theta-unknown} are stable under those transformations made in order to simplify our study.

We define $\mathcal{D}^0(\R^2)$ the set of continuous and bijective functions from $\R^2$ to $\R^2$ with a continuous inverse, taking value 0 at 0.
We define $\mathcal{D}^2(\R^2)$ the set of $\mathcal{C}^2$-diffeomorphisms from $\R^2$ to $\R^2$ taking value 0 at 0. We call such functions (in $\mathcal{D}^0(\R^2)$ or in $\mathcal{D}^2(\R^2)$, according to the section of this paper) deformations.

Note that the assumptions on $X$ and on $\theta$ that we have just listed are not all in force in Section \ref{sectionstatiso}, where we soften the regularity assumptions on $X$ and $\theta$ and we replace the Gaussian hypothesis on $X$ by the assumption of the existence of a second moment.
\section{For which $\theta$ is $X_{\theta}$ isotropic?}
\label{sectionstatiso}
In this section, assumption \textbf{(H)} is \textbf{not} in force. We only assume that $X$ is stationary, isotropic and that it admits a second moment.
We denote by $C_{\theta}$ the covariance function of the deformed field $X_{\theta}$. Because the field $X$ is stationary, for any $(x,y)\in(\R^2)^2$,
\begin{equation}
\label{rthetageneral}
C_{\theta}(x,y)=\Cov(X_{\theta}(x),X_{\theta}(y))=C(\theta(x)-\theta(y)).
\end{equation}

In the following, we exhibit the deformations $\theta$ that leave the field $X_{\theta}$ isotropic, for any stationary and isotropic field $X$. Note that the underlying field $X$ is not fixed. Our approach is analogous to the one in \cite{PS00}, where the objective is, starting with a random field $Y$ with a known covariance function, to find a deformation $\theta$ such that $Y=X\circ\theta$, with $X\,:\,\R^2\rightarrow\R$ a stationary, or stationary and isotropic random field. 

We begin with a short introduction of notations relative to polar representation. We write $\mathbb{T}=\R/2\pi\Z$ the one-dimensional torus and we denote by $S$ the transformation of polar coordinates to cartesian coordinates in the plane deprived of the origin: 
\[S:\; (0,+\infty)\times \mathbb{T}\rightarrow \R^2\backslash\{0\}\quad (r,\varphi)\mapsto(r\cos\varphi,r\sin \varphi).
\]
We define $\mathcal{D}^0\left((0,+\infty)\times \mathbb{T}\right)$ the set of continuous and bijective functions $\hat{\theta}:\,(0,+\infty)\times \mathbb{T}\rightarrow(0,+\infty)\times \mathbb{T}$ with continuous inverses.
For any deformation $\theta\in\mathcal{D}^0\left(\R^2\right)$, we write $\theta_0=\theta_{|\R^2\backslash\{0\}}$, we define the deformation $\hat{\theta}\in\mathcal{D}^0\left((0,+\infty)\times \mathbb{T}\right)$ by $\hat{\theta}=S^{-1}\circ\theta_0 \circ S$ and we denote by $\hat{\theta}_1$ and $\hat{\theta}_2$ its coordinate functions.
\begin{prop}
\label{Propiso}
The application $\mathcal{D}^0\left(\R^2\right)\rightarrow \mathcal{D}^0\left((0,+\infty)\times \mathbb{T}\right)\quad \theta\mapsto \hat{\theta}$ is injective and it is a group morphism, that is to say if $\eta$ and $\theta$ are in the set $\mathcal{D}^0\left(\R^2\right)$ then $\widehat{\eta\circ\theta}=\hat{\eta}\circ\hat{\theta}$. Moreover, the coordinate functions of the composition $\widehat{\eta\circ\theta}$ are
 \[\widehat{\eta\circ\theta}^1=\hat{\eta}^1\circ\hat{\theta} \quad \text{and}\quad \widehat{\eta\circ\theta}^2=\hat{\eta}^2\circ\hat{\theta}.
 \]
 \end{prop}
 {\bf Proof.} The above application is obviously injective and if $\eta$ and $\theta$ belong to $\mathcal{D}^0\left(\R^2\right)$, then
\[(\eta\circ\theta)_0=\eta_0\circ \theta_0=(S\circ\hat{\eta}\circ S^{-1})\circ (S\circ\hat{\theta}\circ S^{-1})=S\circ \hat{\eta}\circ \hat{\theta}\circ S^{-1},
\]
hence we get the homomorphism property. Consequently, for $i\in\{1,2\}$, the coordinate function $\widehat{\eta\circ\theta}^i$ satisfies 
 \[(\widehat{\eta\circ\theta}^1,\widehat{\eta\circ\theta}^2)=\widehat{\eta\circ\theta}=\hat{\eta}\circ
 \hat{\theta}= (\hat{\eta}^1\circ\hat{\theta},\hat{\eta}^2\circ\hat{\theta}).
 \]
 \hfill $\Box$

\begin{defn}
\label{defdefspirale}
A deformation $\theta \in\mathcal{D}^0(\R^2)$ is a spiral deformation if there exist $f:\,(0,+\infty)\rightarrow(0,+\infty)$ strictly increasing and surjective, $g:\,(0,+\infty)\rightarrow\mathbb{T}$ and $\eps\in\{\pm 1\}$ such that $\theta$ satisfies
\begin{equation}
\label{defspiral}
\forall (r,\varphi)\in(0,+\infty)\times \mathbb{T}, \quad \hat{\theta}(r,\varphi)=(f(r), \, g(r)+\eps\varphi).
\end{equation} 
\end{defn}

\begin{rems}
\label{remspir}
Note that the set of spiral deformations forms a group for the operation of composition. The choice of $f$ strictly increasing  is due to the conditions of continuity and inversibility on $\theta$ and to the fact that $\theta(0)=0$. The $2\pi$-periodicity of $\hat{\theta^2}$ entails that coefficient $\eps$ in the angular part of (\ref{defspiral}) is an integer and the inversibility of $\theta$ implies that $\eps$ belongs to $\{\pm 1\}$. If we only consider deformations with positive Jacobian determinants, in accordance with our explanations in Section \ref{sectionnota}, then we can set $\eps=1$. Indeed, the positivity of the Jacobian determinant of $\theta$ is equivalent to the positivity of the one of $\hat{\theta}$ (see Formula (\ref{Jtheta}) in the following by way of justification).
\end{rems}

\begin{ex}
\label{ex-spirale-lin}
\underline{Linear spiral deformations.}
A linear spiral deformation is a deformation with polar representation either $(r,\varphi)\mapsto (\lambda r,\varphi+\alpha)$ or $(r,\varphi)\mapsto (\lambda r,-\varphi+\alpha)$, with $\lambda\neq 0$ and $\alpha\in\mathbb{T}$, that is to say it is of the form $\lambda\rho$, with $\rho\in O(2)$.
\end{ex}

In \cite{C87}, the deformations are restricted to the ones given by symmetric, positive and definite matrices. In that case, the field $X_{\theta}$ is isotropic if and only if the two positive eigenvalues of $\theta$ are equal. In the following theorem, we also determine the deformations preserving isotropy but in the general case.
\begin{thm}
\label{equivspiral}
The deformations in $\mathcal{D}^0(\R^2)$ such that for any stationary and isotropic field $X$, $X_{\theta}$ is isotropic are the spiral deformations.
\end{thm}
{\bf Proof.}
To prove the direct implication, let us assume that a deformation $\theta$ is a spiral deformation with polar representation (\ref{defspiral}) and let $\alpha\in\mathbb{T}$. 
 \begin{equation*}
 \begin{aligned}
 \forall (r,\varphi)\in(0,+\infty)\times \mathbb{T},\quad
 \hat{\theta} \circ \hat{\rho}_{\alpha} &=(f(r),g(r)+\eps(\varphi+\alpha))\\
 &=(f(r),g(r)+\eps\varphi+\eps\alpha)\\
 &=\hat{\rho}_{\eps\alpha}\circ\hat{\theta}.
 \end{aligned}
 \end{equation*}
Therefore, $\theta$ satisfies the following property:
\begin{equation*}
\forall \rho \in SO(2),\quad \exists \rho'\in SO(2)\; /\; \theta\circ\rho=\rho'\circ \theta.
\end{equation*}  
This entails that $X_{\theta}\circ \rho=X\circ\rho'\circ \theta$.  
The isotropy of $X$ implies that $X\circ\rho'$ has the same law as $X$. Consequently, $X_{\theta}\circ \rho$ has the same law as $X_{\theta}$. Thus the isotropy of $X_{\theta}$ is proved.
 
We now turn to the converse implication. Let us assume that for any stationary and isotropic field $X$, the field $X_{\theta}$ is isotropic. Hence its covariance function, given by (\ref{rthetageneral}) is invariant under the action of any rotation: \begin{equation*}
\forall \rho\in SO(2),\;\forall (x,y)\in (\R^2)^2,\quad C_{\theta}(\rho(x),\rho(y))=C_{\theta}(x,y).
\end{equation*}
In particular, if we use the Gaussian covariance function $C(x)=\exp(-\|x\|^2)$, we obtain
\begin{equation}
\label{condcoviso}
\forall \rho\in SO(2),\;\forall (x,y)\in (\R^2)^2,\quad \|\theta(\rho(x))-\theta(\rho(y))\|=\|\theta(x)- \theta(y)\|.
\end{equation}
Taking $y=0$, we deduce from (\ref{condcoviso}) that $\hat{\theta}_1$ is radial. We set for any $\varphi\in\mathbb{T}$ and for any $r>0$, $\hat{\theta}_1(r,\varphi)=f(r)$. Since $\theta$ is bijective, continuous and $\theta(0)=0$, $f$ is necessarily strictly increasing with $\lim_{r\rightarrow 0}f(r)=0$ and $\lim_{r\rightarrow+\infty}f(r)=+\infty$.

To infer the form of $\hat{\theta}_2$, we fix $r>0$ and, for any $\varphi\in \mathbb{T}$, we use the complex representation to write Formula (\ref{condcoviso}) for $x=re^{i\varphi}$, $y=r$ and for any angle $\alpha$ of the rotation $\rho$. Dividing the equality by $f(r)$, we get 
\[|e^{i\hat{\theta}_2(r,\varphi+\alpha)}-e^{i\hat{\theta}_2(r,\alpha)}|=|e^{i\hat{\theta}_2(r,\varphi)}-e^{i\hat{\theta}_2(r,0)}|,
\]
hence
\[|1-e^{i(\hat{\theta}_2(r,\varphi+\alpha)-\hat{\theta}_2(r,\alpha))}|=|1-e^{i(\hat{\theta}_2(r,\varphi)-\hat{\theta}_2(r,0))}|.
\]
Since 1 as well as each exponential term belongs to $\{z\in\mathbb{C}\,/\, |z|=1\}$, a geometric interpretation of the above equality entails that for any $\varphi\in\mathbb{T}$, there exists $\epsilon(r,\varphi,\alpha)\in\{\pm 1\}$ such that 
\begin{equation}
\label{formuleangle}
\hat{\theta}_2(r,\varphi+\alpha)-\hat{\theta}_2(r,\alpha)=\epsilon(r,\varphi,\alpha)\,(\hat{\theta}_2(r,\varphi)-\hat{\theta}_2(r,0)).
\end{equation}
Assuming that there exists $\varphi\neq0$ such that $\hat{\theta}_2(r,\varphi)-\hat{\theta}_2(r,0)=0$, we deduce from (\ref{formuleangle}) that $\hat{\theta}_2(r,\cdot)$ is constant on $\mathbb{T}$, which contradicts the bijectivity of $\theta$. Consequently, for any $\varphi\neq0$,
\[\epsilon(r,\varphi,\alpha)=\frac{\hat{\theta}_2(r,\varphi+\alpha)-\hat{\theta}_2(r,\alpha)}{\hat{\theta}_2(r,\varphi)-\hat{\theta}_2(r,0)}.
\] 
This implies that $\epsilon$ is continuous from $(0,+\infty)\times\mathbb{T}\backslash\{0\}\times \mathbb{T}$ onto $\{\pm 1\}$. A connexity argument applies and implies that $\epsilon$ is constant. We write $\epsilon(r,\varphi,\alpha)=\epsilon\in\{\pm 1\}$.

We fix $r>0$. For any $(\varphi,\alpha)\in\mathbb{T}^2$, we can rewrite (\ref{formuleangle}) 
\[\hat{\theta}_2(r,\varphi+\alpha)=\hat{\theta}_2(r,\alpha)+\epsilon(\hat{\theta}_2(r,\varphi)-\hat{\theta}_2(r,0)).
\]
By differentiating the above equality with respect to $\alpha$, for a fixed $\varphi\in\mathbb{T}$, we deduce that $\partial_{\varphi}\hat{\theta}_2(r,\cdot)$ is constant on $\mathbb{T}$. Therefore, there exists $k(r)\in\{\pm 1\}$ and $g(r)\in\mathbb{T}$ such that 
\[\forall r>0, \;\forall \varphi\in\mathbb{T},\quad 
\hat{\theta}_2(r,\varphi)=k(r)\varphi+g(r).
\]
Note that the reason why $k(r)$ must belong to $\{\pm 1\}$ has already been explained in Remarks \ref{remspir}. Finally, since $\hat{\theta}_2$ is continuous, $k(r)$ is necessarily constant, which concludes the proof of Theorem \ref{equivspiral}.
\vspace{0.3cm}

\begin{rem}
Considering the proof of Theorem \ref{Propiso}, we could state another theorem: the deformations in $\mathcal{D}^0(\R^2)$ such that, if $X$ is a stationary and isotropic field with covariance function $C(x)=\exp(-\|x\|^2)$, $X_{\theta}$ is isotropic, are the spiral deformations. In fact, it would be true for any stationary and isotropic field $X$ with an injective covariance function.
\end{rem}
\section{Euler expectation of an excursion set}
\label{section-Euler}

The Euler characteristic is defined on a large subset of compact sets, the class of basic complexes. There are several ways to define this topological functional (see for instance \cite{L00} and \cite{AT07}). 
However, we are actually only interested in the Euler characteristic of excursion sets, which can be computed thanks to specific formulas. 
From now on, $X$ is a random field assumed to satisfy \textbf{(H)} and $\theta$ is a deformation in $\mathcal{D}^2(\R^2)$. Consequently, even though $X_{\theta}$ is in general not stationary nor isotropic, it is Gaussian and its realizations are almost surely of class $\mathcal{C}^2$. Moreover, if $T$ is a rectangle or a segment in $\R^2$ then the set $\theta(T)$ and its frontier $\delta\theta(T)=\theta(\delta T)$ are compact and piece-wise $\mathcal{C}^2$ manifolds in $\R^2$, of respective dimensions two and one.

We start by introducing the general formula for the Euler characteristic of an excursion set of $X_{\theta}$, above a $d$-dimensional rectangle $T$ and then we show how it adapts to dimensions $d=1$ and $d=2$. Let us first explain why we may study equivalently the stationary and isotropic field $X$ on the transformed set $\theta(T)$ or the non-stationary and anisotropic field $X_{\theta}$ on the set $T$. The deformation $\theta$ is an homeomorphism and it satisfies $A_u(X_{\theta},T)
=\theta^{-1}(A_u(X,\theta(T)))$, therefore the sets $A_u(X_{\theta},T)$ and $A_u(X,\theta(T))$ are homotopic. Since the Euler characteristic is a homotopy invariant (see \cite{L00} Theorem 13.36), the above relation leads to 
\[\chi(A_u(X_{\theta},T))=\chi(A_u(X,\theta(T))).
\]

Consequently, we can focus on $\E[\chi(A_u(X,\theta(T)))]$ that can be computed thanks to \cite{AT07} Theorem 12.4.2. We write $(H_i)_{i\in\N}$ the Hermite polynomials and, for any real $x$, $H_{-1}(x)=\sqrt{2\pi}\Psi(x)\exp(x^2/2)$, where $\Psi$ is the tail probability of a standard Gaussian variable.
\begin{equation}
\label{espe2}
\E[\chi(A_u(X_{\theta},T))]=\E[\chi(A_u(X,\theta(T)))]=\underset{0 \le i \le d}{\sum}\mathcal{L}_i(\theta(T)) \, \rho_i(u),
\end{equation}
\[\text{with}\quad \forall \,0\le i \le d,\;\rho_i(u)=(2\pi)^{-(i+1)/2}H_{i-1}(u)e^{-u^2/2}\]
and with $\mathcal{L}_i(\theta(T))$ the $i^{th}$ Lipschitz-Killing curvature of $\theta(T)$. Thanks to the isotropy assumption on $X$ and to the hypothesis $C''(0)=-I_2$, the Lipschitz-Killing curvatures have a very simple expression (see \cite{AT07} section 12.5: ``Isotropic Fields over Smooth Domains"): 
\begin{equation*}
\begin{aligned}
&\text{if }d=1,\;\mathcal{L}_{1}(\theta(T))=|\theta(T)|_{1},\; \mathcal{L}_0(\theta(T))=\chi(\theta(T))=1,\\
&\text{if }d=2,\; \mathcal{L}_2(\theta(T))=|\theta(T)|_2, \; \mathcal{L}_{1}(\theta(T))=\frac{1}{2}|\partial \theta(T)|_{1},\; \mathcal{L}_0(\theta(T))=\chi(\theta(T))=1.
\end{aligned}
\end{equation*} 
Thus, for any two-dimensional rectangle $T\subset \R^2$, from Formula (\ref{espe2}), we get
\begin{equation}
\label{expectedECdim2}
\E[\chi(A_u(X,\theta(T)))]=e^{-u^2/2}\left(u\frac{|\theta(T)|_2}{(2\pi)^{3/2}}+\frac{|\partial \theta(T)|_{1}}{4\pi }\right)+\Psi(u),
\end{equation}
If $T$ is a segment in $\R^2$ then $\theta(T)$ is a one-dimensional manifold and we apply Formula (\ref{espe2}) with $d=1$:
 \begin{equation}
\label{expectedECdim1}
\E[\chi(A_u(X,\theta(T)))]=e^{-u^2/2}\frac{|\theta(T)|_1}{2\pi}+\Psi(u).
\end{equation}
\sloppy
For our approach in Section \ref{section-theta-unknown}, where we want to identify $\theta$ by considering some well-chosen excursion sets of $X_{\theta}$, it will be easier to limit ourselves to the term of highest index in (\ref{espe2}), called the modified Euler characteristic, and also used in \cite{EL16} and \cite{DEL16}. That is why in the following, our results will also involve the modified Euler characteristic (denoted by $\phi$) of excursion sets.  
The general formula for the expectation of the modified Euler characteristic of an excursion set is 
\begin{equation*}
\E[\phi(A_u(X,\theta(T)))]=\mathcal{L}_d(\theta(T)) \, \rho_d(u),
\end{equation*}
from which we deduce, if $d=2$,
\begin{equation}
\label{expectedMECdim2}
\E[\phi(A_u(X,\theta(T)))]=e^{-u^2/2}u\frac{|\theta(T)|_2}{(2\pi)^{3/2}},
\end{equation}
and if $d=1$,
 \begin{equation}
\label{expectedMECdim1}
\E[\phi(A_u(X,\theta(T)))]=e^{-u^2/2}\frac{|\theta(T)|_1}{2\pi}.
\end{equation}

\begin{rem}[Additivity property]
\label{rem_additive} 
The Euler characteristic is an additive functional, which implies that if $T$ and $T'$ are two regular compact sets in $\R^2$ such that $T\cap T'=\emptyset$ then
\[\mathbb{E}[\chi(A_u(X_{\theta},T\cup T'))]=
\mathbb{E}[\chi(A_u(X_{\theta},T))]
+\mathbb{E}[\chi(A_u(X_{\theta},T'))].\] 
In fact, $\phi$ satisfies the same property because the modified Euler characteristic, as well as the Euler characteristic of an excursion set may be expressed as the alternate sum of numbers of critical points of different types of $X$ in the considered domain (see \cite{AT07} Corollary 9.3.5; the term of highest index in the sum corresponds to the modified Euler characteristic of the excursion set). Even if the two-dimensional sets $T$ and $T'$ have a non-empty but one-dimensional intersection, the additivity property is still satisfied. Indeed, in this case, according to Bulinskaya lemma (\cite{AT07} Lemma 11.2.10), almost surely, $X$ admits no critical points in $T\cap T'$; consequently, the (modified) Euler characteristic of the excursion set of $X$ over $T\cap T'$ is almost surely 0.
\end{rem}

We now state a continuity result on the mean Euler characteristic of excursion sets.  
The proposition hereafter shows that if $T$ is a segment in $\R^2$, the mean Euler characteristic of the excursion set of $X_{\theta}$ above $T$ may be seen as the limit of the mean Euler characteristic of excursion sets of $X_{\theta}$ over a sequence of two-dimensional sets, decreasing in the sense of set inclusion and approaching $T$.

\begin{prop}
\label{PropCVtubes}
Let $T$ be a segment in $\R^2$. Let $u$ be a unit vector orthogonal to $T$ and, for any $\rho>0$, let $T_{\rho}$ be the rectangle $\{t+\delta u,\,t\in T,\,-\rho\le \delta\le \rho\}$. Then, for any random field $X$ satisfying Assumption $\mathbf{(H)}$, as $\rho$ decreases towards 0,
\begin{equation*}
\E[\chi(A_u(X_{\theta},T_{\rho}))]\underset{\rho\rightarrow 0}{\longrightarrow}\E[\chi(A_u(X_{\theta},T))]
\end{equation*}
\end{prop}
{\bf Proof.}
The set $\theta(T)$ is one-dimensional while for any $\rho>0$, $\theta(T_{\rho})$ is two-dimensional. Therefore, according to (\ref{expectedECdim1}) and to (\ref{expectedECdim2}),
\begin{equation*}
\begin{aligned}
\E[\chi(A_u(X,\theta(T))]&=e^{-u^2/2}\frac{|\theta(T)|_1}{2\pi}+\Psi(u),\\
\forall \rho>0,\quad 
\E[\chi(A_u(X,\theta(T{_\rho})))]&=e^{-u^2/2}\left(u\frac{|\theta(T{_\rho})|_2}{(2\pi)^{3/2}}+\frac{|\partial \theta(T{_\rho})|_{1}}{4\pi }\right)+\Psi(u).
\end{aligned}
\end{equation*}

For any sequence $(\rho_n)_{n\in\N}$ of positive terms decreasing towards 0, the sequence of sets $(\theta(T_{\rho_n}))_{n\in\N}$,  decreases to $\cap_{n\in\N}\theta(T_{\rho_n})=\theta(T)$ thus the limit of $|\theta(T_{\rho_n})|_2$ as $n$ tends to infinity is zero.
  
For any $\rho>0$, the frontier of $\theta(T_{\rho})$ is 
\[\partial\theta(T_{\rho})=\theta(\partial T_{\rho})
=\{\theta(t+\rho u),\,t\in T\}\cup
\{\theta(t-\rho u),\,t\in T\}.
\] 
As $\rho$ tends to 0, the one-dimensional measure of each set of this disjoint union tends to $|\theta(T)|_1$; therefore,  
$|\partial \theta(T_{\rho})|_1$ tends to $2|\theta(T)|_1$. This concludes the proof.

\begin{rems}
Proposition \ref{PropCVtubes} could be adapted in various ways. First, we could generalize it to a one-dimensional compact and connected set $T$ satisfying certain regularity assumptions. Besides, the sequence of sets $\{T_{\rho},\,\rho>0\}$ approaching $T$ could be defined differently, for instance as the sequence of $\rho$-tubes around $T$, that is,
\[\forall \rho>0,\,T_{\rho}=\{z\in\R^2\,/\,\text{dist}(T,z)\le \rho\},
\;\text{where }\text{dist}(T,z)=\min_{x\in T}\{\|x-z\|\}.
\]
We should also point out that Proposition \ref{PropCVtubes} is specific to the Euler characteristic and the same result would not stand with the modified Euler characteristic. Indeed, consider for instance $T_N=[a,b]\times [-N^{-1},N^{-1}]$ for $N\in\N\backslash\{0\}$, then Formula (\ref{expectedMECdim2}) yields
\[\E[\phi(A_u(X,\theta(T_N)))]\underset{N\rightarrow 0}{\longrightarrow}0,
\]
whereas, according to Proposition \ref{PropCVtubes},
\[\E[\chi(A_u(X,\theta(T_N)))]\underset{N\rightarrow 0}{\longrightarrow}\exp(-u^2/2)\frac{|\theta([a,b]\times\{0\})|_1}{2\pi}+\Psi(u).
\]
\end{rems} 
To end with, here is an integral formula giving the second moment of the modified Euler characteristic of an excursion set of $X$ over $\theta(T)$, where $T$ is a rectangle. It will be useful in Section \ref{subsec_spiral} when we address estimation matters. It is proved in \cite{DEL16} Proposition 1 for an excursion set over a cube, however it remains true in our case. 
\begin{equation}
\begin{aligned}
\label{2momentMEC}
\Var[\phi(A_u(X_{\theta},T)]&=\Var[\phi(A_u(X,\theta(T))]\\
&=\int_{\R^2}|\theta(T)\cap (\theta(T)-t)|_2(G(u,t)D(t)^{1/2}-h(u)^2)\,dt\\
&\qquad \qquad+|\theta(T)|_2(2\pi)^{-1}g(u),
\end{aligned}
\end{equation}
where 
\begin{equation*}
\begin{aligned}
&G(u,t)=\E[\mathbf{1}_{[u,+\infty)}(X(0))\mathbf{1}_{[u,+\infty)}(X(t))\det(X''(0))\det(X''(t))|X'(0)=X'(t)=0],\\
&D(t)=(2\pi)^{4}\det(I_2-C''(t)^2),\\
&g(u)=\E[ \mathbf{1}_{[u,+\infty)}(X(0))|\det(X''(0))|)],\\
&h(u)= (2\pi)^{-3/2}\,u\,e^{-u^2/2}.
\end{aligned}
\end{equation*}

\section{Notion of $\chi$-isotropic deformation}
\label{section-stariso}
In this section, the underlying field $X$ is fixed and it satisfies Assumption $\textbf{(H)}$. We define $\chi$-isotropic deformations, characterized by an invariance condition of the mean Euler characteristic of some excursion sets of the associated deformed field. We show that the only deformations that satisfy this invariance property are the spiral deformations, that is to say the ones that were proved to preserve isotropy in Section \ref{sectionstatiso}.

\begin{defn}[$\chi$-isotropic deformation]
\label{defchiiso}
A deformation $\theta\in\mathcal{D}^2(\R^2)$ is $\chi$-isotropic if for any rectangle $T$ in $\R^2$, for any $u\in\R$ and for any $\rho \in SO(2)$,
 \begin{equation}
\label{equa_*iso}
\E[\chi(A_u(X_{\theta},\rho(T))]=\E[\chi(A_u(X_{\theta},T)].
\end{equation}
\end{defn}
\begin{rem}
If we consider the modified Euler characteristic instead of the Euler characteristic, we can easily adapt Definition \ref{defchiiso} in the following way: 
a deformation $\theta\in\mathcal{D}^2(\R^2)$ is $\chi$-isotropic if for any rectangle \underline{or segment} $T$ in $\R^2$, for any $u\in\R$ and for any $\rho \in SO(2)$,
 \begin{equation*}
\label{equa_*isoMEC}
\E[\phi(A_u(X_{\theta},\rho(T))]=\E[\phi(A_u(X_{\theta},T)].
\end{equation*}
\end{rem}
\begin{ex}
\label{ex_spiral_chiiso}
Spiral deformations defined in Section \ref{sectionstatiso} are $\chi$-isotropic deformations. Indeed, if a deformation $\theta$ is such that $X_{\theta}$ is isotropic then it satisfies the above definition, because for any $\rho\in SO(2)$, $X_{\theta\circ\rho}$ has the same law as $X_{\theta}$. But according to Theorem \ref{equivspiral}, the deformations that preserve isotropy are exactly the spiral deformations.
\end{ex}

Here comes the main result of Section \ref{section-stariso}.
\begin{thm}
\label{thm_chiiso}
The $\chi$-isotropic deformations are exactly the spiral deformations in $\mathcal{D}^2(\R^2)$.
\end{thm}

Before turning to the proof of Theorem \ref{thm_chiiso}, note that it shows that $\chi$-isotropy is independent on the underlying field $X$, even though the $\chi$-isotropy definition (Definition \ref{defchiiso}) is stated with a fixed one.\\

\noindent
{\bf Proof of Theorem \ref{thm_chiiso}.}
Spiral deformations are $\chi$-isotropic deformations according to Example \ref{ex_spiral_chiiso}; we prove that they are the only $\chi$-isotropic deformations thanks to two lemmas and one result from \cite{BF17}.
The first lemma gives a characterization of $\chi$-isotropic deformations involving invariance properties of the Jacobian matrix under rotations.
\begin{lem}
\label{lem1*iso}
A deformation $\theta\in\mathcal{D}^2(\R^2)$ is $\chi$-isotropic if and only if for any $\rho\in SO(2)$, for any $x\in\R^2$,
\begin{equation}
\label{cond*isot2}
\left\{
\begin{aligned}
(i)&\quad\forall i\in\{1,2\},\;\|J^i_{\theta\circ\rho}(x)\|=\|J^i_{\theta}(x)\|,\\
(ii)&\quad\det(J_{\theta\circ\rho}(x))=\det(J_{\theta}(x)).
\end{aligned}
\right.
\end{equation}
\end{lem}
{\bf Proof}
Let $\theta\in\mathcal{D}^2(\R^2)$ be a $\chi$-isotropic deformation and let $\rho\in SO(2)$. We fix $(s,t)\in\R^2$ and $u\in\R\backslash\{0\}$. Identity (\ref{equa_*iso}) is satisfied for rectangle $T=T(s,t)$, thus Formula (\ref{expectedECdim2}) applied at two different levels $u$ and $u'$ implies that $|\theta\circ\rho(T(s,t))|_2=|\theta(T(s,t))|_2$, whence
\[\int_{[0,s]}\int_{[0,t]} |\det(J_{\theta\circ\rho}(x,y))|\,dx\,dy=
\int_{[0,s]}\int_{[0,t]} |\det(J_{\theta}(x,y)|\,dx\,dy.
\] 
Differentiating twice the above equality with respect to $s$ and to $t$ yields for any $(s,t)\in\R^2$, $|\det(J_{\theta\circ\rho}(s,t))|=|\det(J_{\theta}(s,t))|$, but $|\det(J_{\theta\circ\rho}(s,t))|=
|\det(J_{\theta}(\rho(s,t)))|$ and the Jacobian determinant of $\theta$ has a fixed sign on $\R^2$, hence (\ref{cond*isot2}) Condition \textit{(ii)} is satisfied. Now let us prove Condition \textit{(i)}, for $i=1$ for instance. For any $n\in\mathbb{N}\backslash\{0\}$, according to the definition of $\chi$-isotropy,
\[\E[\chi(A_u(X_{\theta},[0,s]\times[t-n^{-1},t+n^{-1}))]=
\E[\chi(A_u(X_{\theta},\rho([0,s]\times[t-n^{-1},t+n^{-1}])))].
\]
Then we apply Proposition \ref{PropCVtubes} to the set  $[0,s]\times\{t\}$, intersection of the sets $\{[0,s]\times[t-n^{-1},t+n^{-1}],\,n\in\N\}$ (respectively to the set  $\rho([0,s]\times\{t\})$, intersection of the sets $\{\rho([0,s]\times[t-n^{-1},t+n^{-1}]),\,n\in\N\}$). This yields
\[\E[\chi(A_u(X_{\theta},[0,s]\times\{t\}))]=
\E[\chi(A_u(X_{\theta},\rho([0,s]\times\{t\})))],
\]
and, from Formula (\ref{expectedECdim1}),
\[|\theta\circ\rho([0,s]\times\{t\})|_1=|\theta([0,s]\times\{t\})|_1,
\]
which can be written
\[\int_{[0,s]}\|J^1_{\theta\circ\rho}(x,t)\|\,dx=
\int_{[0,s]}\|J^1_{\theta}(x,t)\|\,dx.
\]
Differentiating this integral equality with respect to $s$, we obtain $\|J_{\theta\circ\rho}^1(s,t)\|=\|J_{\theta}^1(s,t)\|$. Similarly, we get $\|J_{\theta\circ\rho}^2(s,t)\|=\|J_{\theta}^2(s,t)\|$. Hence we have proved the direct implication of Lemma \ref{lem1*iso} and we turn to the converse implication.

Let $T$ be a rectangle in $\R^2$. In the first place, there exist $(s,t)\in\R^2$, $\rho_0\in SO(2)$ and a translation by vector $(a,b)\in\R^2$, denoted by $\tau_{a,b}$, such that $T=\rho_0\circ \tau_{a,b} (T(s,t))$. Let $\theta\in\mathcal{D}^2(\R^2)$ satisfying (\ref{cond*isot2}) for any $\rho\in SO(2)$ and for any $x\in\R^2$. Therefore
\begin{equation*}
\begin{aligned}
|\theta\circ\rho(T)|_2&=|\theta\circ\rho\circ \rho_0 (\tau_{a,b}(T(s,t)))|_2\\
&=\int_{[0,s]}\int_{[0,t]}|\det(J_{\theta\circ\rho\circ \rho_0} (a+x,b+y)|\,dx \,dy\\
&=\int_{[0,s]}\int_{[0,t]}|\det(J_{\theta\circ\rho_0} (a+x,b+y)|\,dx \,dy\\
&=|\theta(T)|_2.
\end{aligned}
\end{equation*}
The third equality results from (\ref{cond*isot2}) Condition \textit{(i)}. Now, we express the perimeter length of $\theta\circ\rho(T)$.
\begin{equation*}
\begin{aligned}
|\partial \theta\circ\rho(T)|_1=&|\partial \theta\circ\rho\circ \rho_0(\tau_{a,b}(T))|_1\\
=&\int_{[0,s]} \|J_{\theta\circ\rho\circ \rho_0}^1
(a+x,b)\|\,dx+\int_{[0,s]} \|J_{\theta\circ\rho\circ \rho_0}^1(a+x,b+t)\|\,dx\\
&+\int_{[0,t]} \|J_{\theta\circ\rho\circ \rho_0}^2(a,b+y)\|\,dy+
\int_{[0,t]} \|J_{\theta\circ\rho\circ \rho_0}^2(a+s,b+y)\|\,dy\\
=&\int_{[0,s]} \|J_{\theta\circ \rho_0}^1
(a+x,b)\|\,dx+\int_{[0,s]} \|J_{\theta\circ \rho_0}^1(a+x,b+t)\|\,dx\\
&+\int_{[0,t]} \|J_{\theta\circ \rho_0}^2(a,b+y)\|\,dy+
\int_{[0,t]} \|J_{\theta\circ \rho_0}^2(a+s,b+y)\|\,dy\\
=&|\partial \theta(T)|_1.
\end{aligned}
\end{equation*}
The third equality results from (\ref{cond*isot2}) Condition \textit{(ii)}.
Thanks to Formula (\ref{expectedECdim2}), this proves that $\E[\chi(A_u(X,\theta\circ\rho(T)))]=\E[\chi(A_u(X,\theta(T)))]$. Hence $\theta$ is a $\chi$-isotropic deformation and the proof of Lemma \ref{lem1*iso} is completed.
\hfill
$\Box$
\\~\\
If $M$ and $N$ are two square matrices of size $2\times 2$, the conditions $\|M^i\|=\|N^i\|$ for $i\in\{1,2\}$ and $\det(M)=\det(N)$ are equivalent to the existence of a rotation matrix $\rho\in SO(2)$ such that $M=\rho\,N$. Consequently, we introduce the following equivalence relation on the space of invertible matrices of size $2\times 2$:
\begin{equation*}
\label{equivrel}
M\overset{SO(2)}{\sim}N\quad\Leftrightarrow\quad\exists \rho\in SO(2)\;/\; M=\rho\,N.
\end{equation*}
With this new notation, we can reformulate Lemma \ref{lem1*iso} in the following way: a deformation $\theta\in\mathcal{D}^2(\R^2)$ is $\chi$-isotropic if and only if 
\begin{equation}
\label{carac*isotequiv}
\forall x\in\R^2,\;\forall \rho \in SO(2),\;J_{\theta\circ\rho}(x)\overset{SO(2)}{\sim}J_{\theta}(x).
\end{equation}
Now we are able to state a second lemma that gives another characterization of $\chi$-isotropic deformations involving the polar representation.
\begin{lem}
\label{lem2-*iso}
A deformation $\theta\in\mathcal{D}^2(\R^2)$ is a $\chi$-isotropic deformation if and only if functions
\begin{equation}
\label{syst_EDP}
\left\{
\begin{aligned}
&(r,\varphi)\mapsto(\partial_r \hat{\theta}_1(r,\varphi))^2+(\hat{\theta}_1(r,\varphi)\,\partial_r\hat{\theta}_2(r,\varphi))^2 \\
&(r,\varphi)\mapsto(\partial_{\varphi} \hat{\theta}_1(r,\varphi))^2+(\hat{\theta}_1(r,\varphi)\,\partial_{\varphi}\hat{\theta}_2(r,\varphi))^2\\
&(r,\varphi)\mapsto\hat{\theta}_1(r,\varphi)\,\det(J_{\hat{\theta}}(r,\varphi))\\
\end{aligned}
\right.
\end{equation}
are radial, i.e. if they do not depend on $\varphi$.
\end{lem}
{\bf Proof}
We use the notations introduced at the beginning of Section \ref{sectionstatiso}.
The Jacobian matrix of $S$ at point $(r,\varphi)\in(0,+\infty)\times\mathbb{T}$ is
\[J_S(r,\varphi)=\rho_{\varphi}\begin{pmatrix}
1&0\\0&r
\end{pmatrix}.
\]
Consequently,
\[J_{S^{-1}}(S(r,\varphi))=\left(J_S(r,\varphi)\right)^{-1}=
\begin{pmatrix}
1 &0\\0&r^{-1}
\end{pmatrix}
\rho_{-\varphi}.
\]
Now for any rotation $\rho\in SO(2)$ and for any $(r,\varphi)\in (0,+\infty)\times \mathbb{T}$, we want to express $J_{(\theta\circ\rho)_0}(S(r,\varphi))=J_{\theta_0\circ\rho_0}(S(r,\varphi))$ thanks to $J_{\widehat{\theta\circ\rho}}(r,\varphi)$.

Since $\theta_0=S\circ\hat{\theta}\circ S^{-1}$, we get
\begin{equation}
\label{Jtheta}
J_{\theta_0}(S(r,\varphi))=\rho_{\hat{\theta}_2(r,\varphi)}
\begin{pmatrix}
1&0\\0&\hat{\theta}_1(r,\varphi)
\end{pmatrix}
J_{\hat{\theta}}(r,\varphi)
\begin{pmatrix}
1&0\\0&r^{-1}
\end{pmatrix}
\rho_{-\varphi}.
\end{equation}

 We use the characterization of $\chi$-isotropy given by (\ref{carac*isotequiv}). A deformation $\theta\in\mathcal{D}^2(\R^2)$ is a $\chi$-isotropic deformation if and only if for any $(r,\varphi,\alpha)\in (0,+\infty)\times \mathbb{T}^2$, $J_{\theta_0\circ\rho_{\alpha}}(S(r,\varphi))=J_{\theta_0}(S(r,\varphi+\alpha))\,\rho_{\alpha}$ is equivalent to $J_{\theta_0}(S(r,\varphi))$. Equivalently, for any $(r,\varphi,\alpha)\in (0,+\infty)\times \mathbb{T}^2$, the equivalence relation
\begin{equation*}
\begin{pmatrix}
1&0\\0&\hat{\theta}_1(r,\varphi+\alpha)
\end{pmatrix}
J_{\hat{\theta}}(r,\varphi+\alpha)
\overset{SO(2)}{\sim}
\begin{pmatrix}
1&0\\0&\hat{\theta}_1(r,\varphi)
\end{pmatrix}
J_{\hat{\theta}}(r,\varphi)
\end{equation*}
holds and the above matrices have the same determinant in absolute value and the same norm of columns, which means that functions defined by (\ref{syst_EDP}) do not depend on their second variable.

 \hfill $\Box$
\\~\\
To conclude the proof of Theorem \ref{thm_chiiso}, we refer to \cite{BF17} that uses partial differential equations techniques to prove that deformations satisfying (\ref{syst_EDP}) are spiral deformations.
\hfill $\Box$
\\~\\
Now, we denote by $\mathcal{S}$ the set of spiral deformations in $\mathcal{D}^2(\R^2)$. We also write $\mathcal{I}$ the set of deformations $\theta$ in $\mathcal{D}^2(\R^2)$ such that for any isotropic and stationary field $X$ satisfying $\mathbf{(H)}$, $X_{\theta}$ is isotropic and, for a fixed stationary and isotropic field $X$ satisfying $\mathbf{(H)}$, $\mathcal{I}(X)$ the set of deformations $\theta$ in $\mathcal{D}^2(\R^2)$ such that $X_{\theta}$ is isotropic. Finally, let $\mathcal{X}$ be the set of $\chi$-isotropic deformations. Theorem \ref{thm_chiiso} shows that $\mathcal{X}$ is in fact independent on the underlying field $X$ used in Definition \ref{defchiiso}. The four sets that have just been defined satisfy the following chain of inclusions or equalities:
\[\mathcal{S}=\mathcal{I}\subset\mathcal{I}(X)\subset\mathcal{X}=\mathcal{S}.
\] 
The first and the last equalities come respectively from Theorem \ref{Propiso} and Theorem \ref{thm_chiiso}; the first inclusion is obvious and the second one is a consequence of Example \ref{ex_spiral_chiiso}. As a result, the following corollary holds.
\begin{cor}
\label{Coro}
Let $X$ be a stationary and isotropic random field satisfying Assumption $\mathbf{(H)}$. Then $\mathcal{S}=\mathcal{I}(X)=\mathcal{I}=\mathcal{X}$.
\end{cor}
To conclude, it occurs that the different notions that we have introduced so far to describe the isotropic behaviour of a deterministic deformation are in fact one and correspond to the spiral case.  
\section{Identification of $\theta$ through excursion sets}
\label{section-theta-unknown}
As explained in the introduction of this paper, we consider the case of an unknown deformation $\theta$, which we want to identify using sparse data: the observations of excursion sets of $X_{\theta}$ over well-chosen domains. More precisely, we assume that the mean modified Euler characteristic of some excursion sets of $X_{\theta}$ has been computed and we explain how we can almost uniquely characterize $\theta$. The modified Euler characteristic is more adapted to our method than the Euler characteristic itself. This is due to the dependence of the mean Euler characteristic of an excursion set over a two-dimensional domain on both the perimeter length and the area of the domain, whereas its mean modified version only depends on the area (compare Formulas (\ref{expectedECdim2}) and (\ref{expectedMECdim2})). In the second place, we limit ourselves to spiral deformations and we show that in this case, we can easily estimate $\theta$ thanks to only one realization of the deformed field $X_{\theta}$.

The underlying field $X$ is unknown but it is still assumed to satisfy assumption \textbf{(H)}. The unknown deformation $\theta$ belongs to $\mathcal{D}^2(\R^2)$ and at each point in $\R^2$, its Jacobian determinant is positive. 

\subsection{Identification of $\theta$}
\label{Subsec-charac-theta}
\subsubsection{Case of a linear deformation.}
\label{subsec_linear}
Here comes the simple case of a linear deformation that we use as a first step towards the general case. 
Let us assume that $\theta$ is a linear function and let us write it matricially in a fixed orthonormal basis of $\R^2$: $\theta=\begin{pmatrix}
\theta_{11}&\theta_{12}\\ \theta_{21}&\theta_{22}
\end{pmatrix}$. In this case, we only have to consider the excursion sets over one horizontal segment, one vertical segment and one rectangle (product of two segments): we fix $(s,t)\in(\R\backslash\{0\})^2$, $u\neq 0$ and we assume that we know $\E[\phi(A_u(X_{\theta},[0,s]\times\{0\}))]$, $\E[\phi(A_u(X_{\theta},\{0\}\times[0,t]))]$ and $\E[\phi(A_u(X_{\theta},T(s,t))]$. 
The three real numbers
\begin{equation}
\label{abc}
a=\sqrt{\theta_{11}^2+\theta_{21}^2},\quad
b=\sqrt{\theta_{12}^2+\theta_{22}^2}\quad\text{and}\quad
c=\theta_{11}\theta_{22}-\theta_{21}\theta_{12}
\end{equation}
satisfy
\[
|\theta([0,s]\times \{0\})|_1=|s|a,\quad 
|\theta(\{0\}\times[0,t])|_1=|t|b,\quad \text{and}\quad
|\theta(T(s,t))|_2=|st|c.
\]
Therefore, they are solutions of equations given by Formulas (\ref{expectedECdim1}) and (\ref{expectedECdim2}) and they can be used to write another expression of matrix $\theta$: there exists $(\alpha,\beta)\in\mathbb{T}^2$ such that $\theta=\begin{pmatrix}
a\cos(\alpha)&b\cos(\beta)\\
a\sin(\alpha)&b\sin(\beta)
\end{pmatrix}$.
Let $\delta=\beta-\alpha$ be the angle between the two column vectors; it satisfies $c=ab\sin(\delta)$, whence 
\begin{equation*}
\delta\in\{\delta_0,\delta_1\},\text{ where } \left\{
\begin{aligned}\delta_0&=\arcsin\left(\frac{c}{ab}\right)\in(0,\pi/2]\\\delta_1&=\pi-
\arcsin\left(\frac{c}{ab}\right)\in[\pi/2,\pi).
\end{aligned}
\right.
\end{equation*}
 Consequently, we are able to determine matrix $\theta$ up to an unknown rotation, with two possibilities concerning the angle between its two column vectors: $\theta$ belongs to the set $\mathcal{M}(a,b,c)$ defined by
\begin{equation}
\label{M(abc)}
\mathcal{M}(a,b,c)=\left\{\rho_{\alpha}\begin{pmatrix}
a&\sqrt{b^2-(ca^{-1})^{2}}\\0&ca^{-1}
\end{pmatrix},\;
\rho\begin{pmatrix}
a&-\sqrt{b^2-(ca^{-1})^{2}}\\0&ca^{-1}
\end{pmatrix},\;
\alpha\in \mathbb{T}
\right\}
\end{equation}

If the determinant of $\theta$ was not assumed to be positive, there would be two other possibilities, up to a rotation, because $\delta$ could take four possible values. Note that according to Example \ref{ex-spirale-lin}, $X_{\theta}$ is isotropic in the case where $a=b=\sqrt{c}$, which implies $\delta=\pi/2$.

Of course, because of the isotropy of $X$, we obtain $\theta$ up to post-composition with an unknown rotation. Our method is based on the mean Euler characteristic of excursion sets of $X_{\theta}$ over some sets, which only depends on $\theta$ through the perimeter and area of the set's image by $\theta$. Consequently, we can not differentiate between two deformations that transform any set into sets with the same perimeter and the same area.

 We summarize our approach in the following method.

\begin{met}
\label{proplinearcase}
Let $\theta=\begin{pmatrix}
\theta_{11}&\theta_{12}\\\theta_{21}&\theta_{22}
\end{pmatrix}$ be an unknown linear deformation with positive determinant. For a fixed $(s,t)\in\left(\R\backslash\{0\}\right)^2$, for a fixed $u\in\R\backslash\{0\}$, we assume that $\E[\phi(A_{u}(X_{\theta},T))]$ is known for $T$ of the form $[0,s]\times\{0\}$, $\{0\}\times[0,t]$ and $[0,s]\times[0,t]$. Then $a$, $b$ and $c$ given by (\ref{abc}) are computable thanks to Formulas (\ref{expectedECdim1}) and (\ref{expectedECdim2}) and $\theta$ belongs to the set $\mathcal{M}(a,b,c)$ defined by (\ref{M(abc)}).
\end{met}
\subsubsection{General method.}
\label{subsec_generalapproach}
We refer to the appendix of \cite{AS08} for a precise definition of the complex dilatation and for the statement of the mapping theorem that formulates a characterization of a deformation up to a conformal mapping through its complex dilatation.
To be able to apply it, we add an hypothesis on $\theta$: from now on, we assume that $\theta$ has uniformly bounded distortion, that is to say the ratio of 
\[\displaystyle
\underset{x\rightarrow x_{0}}{\limsup}\frac{|\theta(x)-\theta(x_0)|}{|x-x_0|}\quad\text{to}\quad
\underset{x\rightarrow x_{0}}{\liminf}\frac{|\theta(x)-\theta(x_0)|}{|x-x_0|}
\]
is uniformly bounded for $x_0\in\R^2$.

We fix $u\neq 0$, $S>0$ and we assume that $\E[\phi(A_u(X_{\theta},[0,s]\times\{t\}))]$, $\E[\phi(A_u(X_{\theta},\{s\}\times[0,t]))]$ and $\E[\phi(A_u(X_{\theta},T(s,t))]$ are known for any $(s,t)\in [-S,S]^2$. Then for any $(s,t)\in [-S,S]^2$, we can deduce $|\theta([0,s]\times\{t\})|_1$ and $|\theta(\{s\}\times[0,t])|_1$ from Formula (\ref{expectedMECdim1}) by simply solving a linear system.
Besides
\begin{equation*}
\begin{aligned}
|\theta([0,s]\times\{t\})|_1&= \int_{[0,s]}\|J_{\theta}^{1}(x,t)\|dx =\int_{[0,s]} \sqrt{\partial_x\theta_1(x,t)^2+\partial_x \theta_2(x,t)^2}\, dx,\\
|\theta(\{s\}\times[0,t])|_1&= \int_{[0,t]}\|J_{\theta}^{2}(s,y)\|dy =\int_{[0,t]} \sqrt{\partial_y\theta_1(s,y)^2+\partial_y \theta_2(s,y)^2}\, dy. 
\end{aligned}
\end{equation*}

The first-order partial derivatives of $\theta$ are continuous. By differentiating the functions $s\mapsto |\theta([0,s]\times\{t\})|_1$ and $t\mapsto |\theta(\{s\}\times[0,t])|_1$, we access to functions $s\mapsto \|J_{\theta}^{1}(s,t)\|$ and $t\mapsto \|J_{\theta}^{2}(s,t)\|$ on segment $[-S,S]$.

Now considering the rectangle domains $\{T(s,t),\,(s,t)\in[-S,S]^2\}$, we assume that $\E[\phi(A_u(X,\theta(T(s,t))))]$ is known. Since $u\neq0$, we can compute $|T(s,t)|_2$ thanks to Formula (\ref{expectedMECdim2}).
 Then, by differentiating twice the function
\[(s,t)\mapsto |\theta(T(s,t))|_2=
\int_{[0,s]}\int_{[0,t]} |\det(J_{\theta}(x,y))|\,dx\,dy,
\]
with respect to $s$ and to $t$ on the square $[-S,S]^2$, we obtain function $(s,t)\mapsto|\det(J_{\theta}(s,t))|$ on the same square.

Now, we fix $x\in[-S,S]^2$, we write $J_{\theta}(x)=\begin{pmatrix}
\theta_{11}&\theta_{12}\\ \theta_{21}&\theta_{22}
\end{pmatrix}$ and we use the same notations $a$, $b$ and $c$ defined by (\ref{abc}) as in the linear case, although they now depend on $x$. The explanations given in Section \ref{subsec_linear} apply here and consequently, $J_{\theta}(x)$ belongs to $\mathcal{M}(a,b,c)$.  
Moreover, let us express the complex dilatation $\mu$, given by 
\[\mu=\frac{\partial_{\bar{z}} \theta}{\partial_{z}\theta},
\]
\begin{equation*}
\text{where }
\left\{
\begin{aligned}
&\partial_{z} \theta=\frac{1}{2}(\partial_s \theta_1+\partial_t \theta_2)+\frac{i}{2}(\partial_s \theta_2-\partial_t \theta_1)\\
&\partial_{\bar{z}} \theta=\frac{1}{2}(\partial_s \theta_1-\partial_t \theta_2)+\frac{i}{2}(\partial_s \theta_2+\partial_t \theta_1).
\end{aligned}
\right.
\end{equation*}
At point $x$, a short computation shows that $\mu(x)$ takes two possible values in the set $\mathcal{C}(a,b,c)$ defined by 
\begin{equation}
\label{Cabc}\mathcal{C}(a,b,c)=\left\{\frac{1}{a^2+b^2+2c}(a^2-b^2
\pm 2 i\sqrt{a^2b^2-c^2})\right\}.
\end{equation}
The general method is summarized below.
\begin{met}
\label{propgeneralcase}
\sloppy
Let $\theta\in\mathcal{D}^2(\R^2)$ a deformation with a positive Jacobian on $\R^2$. 
Let $S>0$ and let $u\in\R\backslash\{0\}$ be fixed. Assuming that for any $x=(s,t)\in[-S,S]^2$, for any \mbox{$T\in\left\{[0,s]\times\{t\},\,\{s\}\times[0,t],\,[0,s]\times[0,t]\right\}$}, we know $\E[\phi(A_{u}(X_{\theta},T))]$, we may compute $a=\|J_{\theta}^1(x)\|$, $b=\|J_{\theta}^1(x)\|$ and $c=\det(J_{\theta}(x)$. Consequently, for each $x\in[-S,S]^2$, the Jacobian matrix at point $x$, $J_{\theta}(x)$ belongs to $\mathcal{M}(a,b,c)$ defined by (\ref{M(abc)}) and the complex dilatation at point $x$, $\mu(x)$ belongs to $\mathcal{C}(a,b,c)$ defined by (\ref{Cabc}).
 \end{met}
 \begin{rem}[Numerical approach]
 \sloppy
In practise, we can only have at our disposal a finite amount of data. Let $\sigma$ be a partition of $[-S,S]$. If we know $\left\{\E[\phi(A_{u}(X_{\theta},T))],\;T\in\left\{[0,s]\times\{t\},\,\{s\}\times[0,t],\,[0,s]\times[0,t]\right\},\;(s,t)\in\sigma^2\right\}$, numerical approaches such as Runge-Kutta methods allow to compute approximate values of $\|J_{\theta}^{1}(s,t)\|$, $\|J_{\theta}^2(s,t)\|$ and $\det(J_{\theta}(s,t))$ for any $(s,t)\in\sigma^2$ and the approximate values for $J_{\theta}(s,t)$ and $\mu(s,t)$.
\end{rem}
\subsubsection{Case of a tensorial deformation}
\label{subsec_tensorial}
We now study the particular case of tensorial deformations, where we can completely identify $\theta$ if we make an assumption of monotonicity on its coordinate functions. Let $\theta(s,t)=(\theta_1(s), \theta_2(t))$. Our hypotheses on $\theta$ mean that for $i\in\{1,2\}$, $\theta_i\,:\, \R\rightarrow \R$ satisfies $\theta_i(0)=0$, $\theta_i$ is a bijective function of class $\mathcal{C}^2$ and therefore it is monotonous. Note that $\theta$ transforms a rectangle $[s,v]\times[t,w]$ into another rectangle $\theta_1([s,v])\times \theta_2([t,w])$. 

Let $s\in\R\backslash\{0\}$. We deduce from Formula (\ref{expectedMECdim1}) that
\begin{equation*}
\left\{
\begin{aligned}
\E[\phi(A_u(X_{\theta},[0,s]\times\{0\})]-\Psi(u)&=\frac{e^{-u^2/2}}{2\pi}&\int_{0}^s|\theta_1'(x)|\,dx\\
\E[\phi(A_u(X_{\theta},\{0\}\times[0,s])]-\Psi(u)&=\frac{e^{-u^2/2}}{2\pi}&\int_{0}^s|\theta_2'(x)|\,dx
\end{aligned}
\right.
\end{equation*}
and consequently, we can state the following method.
 \begin{met}
 \label{proptensorialcase}
Let $(s,t)\mapsto\theta(s,t)=(\theta_1(s),\theta_2(t))\in\mathcal{D}^2(\R^2)$ be a tensorial deformation. We fix $S>0$ and $u\in\R$. We assume that for any real number $s\in[-S,S]\backslash \{0\}$ and for $T\in\left\{[0,s]\times\{0\},\,\{0\}\times[0,s]\right\}$, we know $\E[\phi(A_u(X_{\theta},T))]$. Then we determine functions $s\mapsto|\theta_1'(s)|$ and $s \mapsto|\theta_2'(s)|$ on $[-S,S]$ thanks to Formula (\ref{expectedMECdim1}). If the sign of each coordinate function is known then $\theta$ is completely determined.
\end{met}
\begin{ex} Let $(\alpha,\beta)\in(\R\backslash\{0\})^2$, let $\theta$ be defined on $[0,1]^2$ by $\theta(s,t)=(s^{\alpha},t^{\beta})$ and let $\sigma$ be a partition of $(0,1]$. To identify $\theta_1$, we follow the above method adapted to a numerical approach; thus we obtain approximate values for $\{|\theta_1'(s)|,\,s\in \sigma\}$. Constant values correspond to the case of $\alpha=1$. Otherwise, we have $|\theta_1'(s)|=\left|\alpha\right|s^{\alpha-1}$, therefore coefficient $\alpha$ can be computed through a regression method as the slope of the line representing $\log(|\theta_1'(s)|)=\log\left(\left|\alpha\right|\right)+(\alpha-1)\log(s)$
as a function of $\log(s)$ on $(0,1]$. The same method holds to get coefficient $\beta$.
\end{ex}
\begin{rem}
The three methods \ref{proplinearcase}, \ref{propgeneralcase} and \ref{proptensorialcase} can be easily adapted if the modified Euler characteristic $\phi$ is replaced by the Euler characteristic itself $\chi$.
\end{rem}
\subsection{Estimation in the spiral case}
 \label{subsec_spiral}
We have assumed all along the first part of this section that $\mathbb{E}[\phi(A_u(X_{\theta},T))]$ was known for some basic domains $T$, but we have not yet discussed estimation matters. Without any hypothesis on $\theta$, this expectation seems uneasy to estimate from one single realization of $X_{\theta}$, for the deformed field is non-stationary, except in the linear case. Yet it is possible in the spiral case thanks to the isotropy of the deformed field. More precisely, let $\theta\in\mathcal{D}^2(\R^2)$ be a spiral deformation; we show in the following how to estimate $\|J_{\theta}^1(x)\|$, $\|J_{\theta}^2(x)\|$ and $\det(J_{\theta}(x))$ at each point $x$ in a chosen domain. Then the end of Method \ref{propgeneralcase} applies to identify $\theta$.

Let $x\in\R^2\backslash\{0\}$, let $(r_0,\varphi_0)$ be its polar coordinates and for $N\in\N\backslash\{0\}$, let $T_N^0=\{(r,\varphi)\in (0,+\infty)\times \mathbb{T}\;/\;r_0\leq r\leq r_0+N^{-1},\;\varphi_0\le \varphi\le \varphi_0+2\pi N^{-1}\}$.
 
 For any $k\in\{0,\cdots, N-1\}$, we write $T_N^k=\rho_{2k\pi/N}(T_0^N)$. We fix $u\neq 0$ and we define 
\begin{equation*}
Z_N=N^{-1}\sum_{k=0}^{N-1}\phi(A_u(X_{\theta},T_N^k))=
N^{-1}\sum_{k=0}^{N-1}\phi(A_u(X,\theta(T_N^k))),
\end{equation*}
where $\phi$ is the modified Euler characteristic. Remember that $\phi$ satisfies an additivity property (see Remark \ref{rem_additive}). Thus,
\begin{equation*}
\begin{aligned}
Z_N&=N^{-1}\phi\left(A_u\left(X,\underset{k=0}{\overset{N-1}{\cup}}\theta(T_N^k)\right)\right)\\&=N^{-1}\phi(A_u(X,\theta(U_N))),\quad\text{where}\quad
U_N=\underset{k=0}{\overset{N-1}{\cup}}T_N^k.
\end{aligned}
\end{equation*}
We derive the asymptotic behaviour of the expectation and the variance of $Z_N$ from the $\chi$-isotropy property satisfied by $\theta$ according to Theorem \ref{thm_chiiso}.
\begin{prop}
\label{Prop-esp-var-ZN}
There exist constants $a\neq 0$ and $c>0$ (depending only on $u$) and $n\in\N\backslash\{0\}$ such that
\[\mathbb{E}[Z_N]\underset{N\rightarrow+\infty}{\sim}a\,|\det(J_{\theta}(x))||T_N^0|_2\]
and for $N\ge n$, 
\[\Var[Z_N]\le c\,\frac{|\det(J_{\theta}(x))|\,|T_N^0|_2}{N}.
\]
\end{prop}
 \textbf{Proof}
Let $N\in\N\backslash\{0\}$. The $\chi$-isotropy of $\theta$ implies that for any $k\in\{0,\cdots,N-1\}$, $\mathbb{E}[\phi(A_u(X,\theta(T_N^k)))]=\mathbb{E}[\phi(A_u(X,\theta(T_N^0)))]$, according to Definition \ref{defchiiso}. Therefore, the expectation of $Z_N$ is
\[\mathbb{E}[Z_N]=\mathbb{E}[\phi(A_u(X,\theta(T_N^0)))]=\frac{u\,e^{-u^2/2}}{(2\pi)^{3/2}}|\theta(T_N^0)|_2.
1\]
We study the asymptotic behaviour of this sequence.
 \begin{equation*}
 \begin{aligned}
 &\left||\theta(T_N^0)|_2-
 |\det(J_{\theta}(x))||T_N^0|_2\right|\\
 &\le\int_{r_0}^{r_0+N^{-1}}
 \int_{\varphi_0}^{\varphi_0+2\pi N^{-1}}
 \left|\,|\det (J_{\theta}(S(r,\varphi)))|
 -|\det (J_{\theta}(S(r_0,\varphi_0)))| \,\right|r
 \,dr\,d\varphi \\
 &\le \underset{\underset{\varphi_0\le\varphi\le\varphi_0+2\pi N^{-1}}{r_0\le r\le r_0+ N^{-1}}}{\sup}
\left|\,|\det(J_{\theta}(S(r,\varphi)))|-|\det(J_{\theta}(S(r_0,\varphi_0)))|\,\right|\,|T_N^0|_2,
 \end{aligned}
 \end{equation*}
with $\underset{\underset{\varphi_0\le\varphi\le\varphi_0+2\pi N^{-1}}{r_0\le r\le r_0+ N^{-1}}}{\sup}
\left|\,|\det(J_{\theta}(S(r,\varphi)))|-|\det(J_{\theta}(S(r_0,\varphi_0)))|\,\right|\underset{N\rightarrow+\infty}{\longrightarrow}0$. Consequently, 
\begin{equation}
\label{asymptAireImage}
|\theta(T_N^0)|_2\underset{N\rightarrow+\infty}{\sim}|\det(J_{\theta}(x))||T_N^0|_2
\end{equation}
and the result about the asymptotic expectation holds.

Now we use Formula (\ref{2momentMEC}) (with its notations) to get an integration expression of the variance of $Z_N=N^{-1}\phi(A_u(X,\theta(U_N)))$ and an asymptotic upper-bound:
\begin{align*}
\label{majvarsomcour}
\Var[\phi(A_u(X,\theta(U_N)))]&=
\int_{\R^2}|\theta(U_N)\cap (\theta(U_N)-t)|_2(G(u,t)D(t)^{1/2}-h(u)^2)\,dt\\
&\qquad \qquad \qquad+|\theta(U_N)|_2(2\pi)^{-1}g(u)\\
&\le |\theta(U_N)|_2\left(
\int_{\R^2}(G(u,t)D(t)^{1/2}-h(u)^2)\,dt
+(2\pi)^{-1}g(u)\right)\\
&\le c\,|\theta(U_N)|_2=c\,N|\theta(T_N^0)|_2,
\end{align*}
where $c>0$. (Note that the integral domain is in fact the compact $\{t-t',\;(t,t')\in U_N^2\}$.) 
Consequently, using (\ref{asymptAireImage}), we get
\begin{align*}
\Var[Z_N]\leq c\, N^{-1}\,|\theta(T_N^0)|_2
\underset{N\rightarrow+\infty}{\sim}c\,N^{-1}\,|\det(J_{\theta}(x))|\,|T_N^0|_2.
\end{align*}
This concludes the proof.

Proposition \ref{Prop-esp-var-ZN} shows that, asymptotically, the variance of $Z_N$ is negligible with respect to its expectation. Practically, we could obtain $|\det(J_{\theta}(x))|$ through a regression method since, up to a constant, it is the coefficient of the linear relation linking asymptotically $|\theta(T^0_N)|_2$ and $|T_N^0|_2$. Constant $a$ is totally explicit and constant $c$ may be numerically computed.

We can adopt the same approach to get an estimation of $\|J_{\theta}^i(x)\|$, for $i\in\{1,2\}$. We will only state the asymptotic result (for $i=1$) because the proof is very similar to the one of Proposition \ref{Prop-esp-var-ZN}. Let $x=(x_1,x_2)\in\R^2$ and $S^0_N=[x_1,x_1+N^{-1}]\times\{x_2\}$. For $N\in\N$, for any $k\in\{0,\cdots, N-1\}$, we write $S_N^k=\rho_{2k\pi/N}(S^0_N)$ and we define 
\[Y_N=N^{-1}\sum_{k=0}^{N-1}\phi(A_u(X_{\theta},S_N^k)).
\]
\begin{prop}
\label{Prop-esp-var-YN}
There exist constants $d\neq 0$ and $k>0$ (depending only on $u$) and $n\in\N\backslash\{0\}$ such that 
\[\mathbb{E}[Y_N]\underset{N\rightarrow+\infty}{\sim}d\,\|J_{\theta}^1(x)\||S_N^0|_1\]
and for $N\ge n$, 
\[\Var[Y_N]\le k\,\frac{\|J_{\theta}^1(x)\|\,|S_N^0|_1}{N}.
\]
\end{prop}
Estimates of $|\det(J_{\theta})|$, $\|J_{\theta}^1(x)\|$ and $\|J_{\theta}^2(x)\|$ bring a nearly complete characterization of the Jacobian matrix of $\theta$, as explained in Section \ref{Subsec-charac-theta}.
\section*{Acknowledgements}
The author is grateful to her advisor A. Estrade for her valuable help and advice, and to M. Briant for his interest in a collaboration (presented in \cite{BF17}) that substancially helped to produce Theorem \ref{thm_chiiso}.


\end{document}